\documentstyle{article}
\title{Codimension two nonsingular subvarieties of quadrics:\\
scrolls and classification in  degree $d\leq 10$}
\author{Mark Andrea A.  de Cataldo}
\date{}

\newtheorem{tm}{Theorem}[subsection]
\newtheorem{lm}[tm]{Lemma}
\newtheorem{pr}[tm]{Proposition}
\newtheorem{rmk}[tm]{Remark}
\newtheorem{cor}[tm]{Corollary}

\newtheorem{defin}[tm]{Definition}
\newtheorem{?}[tm]{Question}
\newtheorem{evidence}[tm]{}

\newtheorem{fact}[tm]{Fact}

\font\tenmsb=msbm10
\font\sevenmsb=msbm7
\font\fivemsb=msbm5

\newfam\msbfam
\textfont\msbfam=\tenmsb
\scriptfont\msbfam=\sevenmsb
\scriptscriptfont\msbfam=\fivemsb
\def\Bbb#1{{\fam\msbfam #1}}

\font\teneufm=eufm10
\font\seveneufm=eufm7
\font\fiveeufm=eufm5
\newfam\eufmfam
\textfont\eufmfam=\teneufm
\scriptfont\eufmfam=\seveneufm
\scriptscriptfont\eufmfam=\fiveeufm
\def\frak#1{{\fam\eufmfam\relax#1}}

\newcommand\ci{\cite}
\newcommand\G{\Gamma}
\newcommand\s{\sigma}

\newcommand\zed{{\Bbb Z}}
\newcommand\pn[1]{{\Bbb P}^{#1}}

\newcommand{\Q}[1]{{\cal Q}^{#1}}
\newcommand{\Il}[3]{{\cal I}_{{#1},{#2}}({#3})}
\newcommand{\I}[2]{{\cal I}_{{#1},{#2}}}

\newcommand\blacksquare{{\hspace*{\fill} $\Box$}} 
\newcommand\odix[1]{ {\cal O}_{#1} }
\newcommand\odixl[2]{ {\cal O}_{#1}({#2}) }

\newcommand\nb[2]{ {\cal N}_{ {#1}, {#2} } }

\begin{document}
\maketitle

\begin{abstract}
Let $X$ be  a codimension two  nonsingular subvariety of a nonsingular 
quadric
$\Q{n}$ of dimension $n\geq 5$. We classify such subvarieties
when they  are   scrolls. We also classify them when the degree
$d\leq 10$. Both results were known when $n=4$. 
\end{abstract}

\section{INTRODUCTION}
\label{intr}
The paper \ci{ottp5} completes the classification of scrolls as 
codimension 
two subvarieties
of projective space $\pn{n}$. Ottaviani's  proof consists of three parts.
First 
the sectional genus $g$ is exhibited  as  a function of 
the degree $d$ of the scroll.  The degree $d$ is then bounded from 
above  by the use of Castelnuovo-type bounds for $g$.
The final step consists of the construction of  varieties with prescribed
low invariants which  had been accomplished by several authors.

In this paper we classify scrolls as codimension two subvarieties of
$\Q{n}$; see Theorem \ref{maintm}. The analysis is quite similar to 
the one of \ci{ottp5} with the 
following three differences.
The first one is that there are fourfolds scrolls on  $\Q{6}$.
The second difficulty is that  the method
for  bounding  the degree of  scrolls over surfaces on $\Q{5}$ of 
\ci{ottp5}  is  not
sufficient; we go around the problem using 
lemmata \ref{d=6...42} and \ref{d=6812}.
Lastly, once we obtain a maximal list of invariants we must construct
all the scrolls in question.  This is essentially 
the problem of constructing
varieties of low degree  and codimension two on $\Q{n}$.  We build on 
the results of \ci{a-s}
and \ci{gross} and obtain 
Theorem \ref{classificationd<12}, i.e.   the complete classification
in degree $d\leq 10$ and $n\geq 5$. This result highlights the role
that some special vector bundles on quadrics play in  the construction of 
subvarieties
of quadrics.
As a by-pass result of this classification in low degree we are able to 
  construct  all
 scrolls,  except for one case: 
when the  degree $d=12$ and  the base is a minimal $K3$
surface. We construct  an  unirational  family
of these scrolls; see Theorem \ref{esempium}. We do not know whether or 
not 
this   is the only one. 

The paper is organized as follows.
Section \ref{prel} contains preliminary results. 
Section  \ref{d<=10} contains Theorem \ref{classificationd<12},
Section \ref{secpbundles} contains the main result of this paper,
Theorem \ref{maintm}, the proof of which guides the reader through
the rest of the section.

\smallskip \noindent
{\bf Notation and conventions.}
Our basic reference is [Ha].
We work over any algebraically closed field of characteristic zero.
A quadric $\Q{n}$,
here, is a nonsingular hypersurface of degree two in the  projective space
$\pn{n+1}$.
Little or no distinction is made between line bundles, associated 
sheaves of sections
and Cartier divisors.
$\lfloor t\rfloor $
denotes the biggest integer smaller than or equal to 
$t$. $\sim_n$ denotes the numerical equivalence of divisors on a surface.
$\odixl{\Q{n}}{1}$ denotes the sheaf
$\odixl{\pn{n+1}}{1}_{|\Q{n}}$. If $F$ is a coherent sheaf on $\Q{n}$ and
$l$ an integer, 
then $F(l)$ denotes the sheaf
$F\otimes \odixl{\Q{n}}{l}$.

\smallskip
\noindent
{\bf Acknowledgments.} This paper is an expanded
and completed version of parts of our dissertation.
It is a pleasure to thank our Ph.D. advisor
A.J. Sommese, who has suggested to us that we  study threefolds
on $\Q{5}$. We thank the C.N.R. of the Italian Government
and The University of Notre Dame for partial support.
We wish to thank  K. Chandler and J. Migliore
for valuable discussions concerning Proposition
\ref{d=6812} and E. Arrondo and G. Ottaviani  for  useful correspondences.

\section{PRELIMINARY MATERIAL}
\label{prel}

In this section we collect the various results that will be necessary
in sections \ref{d<=10} and \ref{secpbundles}.
In this section $X$ is a codimension two, nonsingular
subvariety of  $\Q{n}$,   $d$ its degree, 
$\nb{X}{\Q{n}}$ its normal bundle, $n_i$ the  $i^{th}$ Chern class
of $\nb{X}{\Q{n}}$,  $\iota:X\hookrightarrow \Q{n}$ the embedding,  
$L$ the restriction of
$\odixl{\Q{n}}{1}$ to $X$, $K_X$ the canonical dualizing sheaf of $X$.

\subsection{Miscellanea}
\label{miscellanea}

The cohomology ring of a  nonsingular quadric of any dimension is
described in \ci{h-p}. Let $h$ be the class of 
any hyperplane section of $\Q{n}$.

We consider the odd dimensional case first: $\Q{2n+1}$.
One can describe $H^*(\Q{2n+1},$ $\zed)$ as follows.
Let
$\Lambda$
be the class of an $n$-dimensional linear space in 
$\Q{2n+1}.$
 The relevant information is, denoting the cup product by
  ``$\cdot$":

\smallskip
\noindent
 $H^{2i+1}(\Q{2n+1},{\Bbb Z})=\{0\}$,  $\forall i$;\,
$H^{2i}(\Q{2n+1},{\Bbb Z})=$ $\{0\}$, for $i>2n+1$;\,
$H^{2i}(\Q{2n+1},{\Bbb Z})=$ ${\Bbb Z}[h^i]$, 
$i=0,\ldots, n$;\, 
$H^{2(n+j)}(\Q{2n+1},{\Bbb Z})$ $={\Bbb Z}[\Lambda\cdot h^{j-1}]$,
$j=1,\ldots,n+1$;\, 
$h^{n+1}=2\Lambda$, $h^{2n+1}=2$.

\smallskip
As to the even dimensional case, we  denote by
$\Lambda_1$, $\Lambda_2$ the classes of two members of the two rulings of
$\Q{2n}$ in $n$-dimensional linear spaces.
One has:

\smallskip
\noindent
 $H^{2i+1}(\Q{2n},{\Bbb Z})=\{0\}$,  $\forall i$;\, 
$H^{2i}(\Q{2n},{\Bbb Z})=$ $\{0\}$, for $i>2n$;\,
$H^{2i}(\Q{2n},{\Bbb Z})=$ ${\Bbb Z}[h^i]$, 
$i=0,\ldots, n-1$;\, 
$H^{2i}(\Q{2n},{\Bbb Z})=$ $\zed [\Lambda_1]\bigoplus \zed[\Lambda_2]$;\, 
$H^{2(n+j)}(\Q{2n+1},{\Bbb Z})$ $={\Bbb Z}[\Lambda_1\cdot h^{j-1}]$ 
$={\Bbb Z}[\Lambda_2\cdot h^{j-1}]$,
$j=1,\ldots,n$;\, 
$h^{n}=\Lambda_1+\Lambda_2$;\, $h^{2n}=2$;\, $[\Lambda_i] \cdot [\Lambda_j]=
 \delta_{ij}$,
where $\delta_{ij}$ is the Kronecker symbol.

\begin{rmk}
\label{rmkdiseven}
{\rm 
   The above description
of the cohomology ring of $\Q{n}$ implies that, for $n\geq 5$, $d$ 
is an even integer.
}
\end{rmk}

 Mumford's {\em self intersection formula}  
(cf. \ci{fu}, page 103) gives, for $n\geq 5$:
\begin{equation}
\label{muself}
n_2=\frac{1}{2}d L^2.       
\end{equation}

Consider the twisted ideal sheaves
$\Il{X}{\Q{n}}{l}:=$ $ \I{X}{\Q{n}}\otimes \odixl{\Q{n}}{l}$. 
 We write the total Chern class of these sheaves as
$1+\sum_{i=1}^{n}\gamma_i h^i$. The following is a standard  a 
consequence of
\ci{fu}, Theorem and Lemma 15.3.

\begin{lm}
\label{lmccis}
Let $X$ and $\Il{X}{\Q{n}}{l}$ be as above, with $l$ fixed.
Assume that $n\geq 5$.
Then one has the following relations concerning the Chern classes
of \,
$\Il{X}{\Q{5}}{l}:$ 
$$
\gamma_1  =  l; \qquad \qquad  
\gamma_i  =  \frac{1}{2}  (K_X+(5-l)L)^{i-2}\cdot L^{n-i}, \quad \forall
i=2,\ldots, n.    
$$
\end{lm}

 We now 
make explicit the  {\em Double Point Formul\ae \ } for the embedding
$\iota$.
 The proof is a standard consequence of \ci{fu} Theorem 9.3, once we use
(\ref{muself}) and the fact that $n_i=0$ for $i\geq 3$.
Denote by $c_i$
the Chern classes of the tangent bundle of $X$.

\begin{lm}
\label{lmDPF}
Let $\iota:X\hookrightarrow \Q{n}$ be as above with
  $n\geq 5$. 
 Then one has the 
following relations in the Chow ring of $X$:

\begin{equation}
\frac{1}{2} d L^2=\frac{1}{2}(n^2 -n +2)L^2 -n c_1\cdot L + c_1^2 - c_2;
\label{deg2dpf5}
\end{equation}
\begin{equation}
c_3= \frac{1}{6}(n^3 -3n^2+8n-12)L^3 +\frac{1}{2}(-n^2 +n -2)c_1 L^2+
n(c_1^2-c_2)L +2c_1c_2-c_1^3 ;
\label{deg3dpf5}
\end{equation}
\begin{equation}
\label{deg4dpf}
c_4= 22L^4 - 24L^3c_1 + 16L^2(c_1^2-c_2) + 12Lc_1c_2 - 6Lc_1^3 -6Lc_3
+2c_1c_3 + (c_1^2-c_2)^2 -c_1^2c_2.
\end{equation}
\end{lm}

For $n=5$ we have:

\begin{equation}
\label{KL2}
K_X\cdot L^2=2(g-1)-2d,
\end{equation}
\begin{equation}
\label{K2L}
K_X^2\cdot L=\frac{1}{4}d^2 +\frac{3}{2}d -8(g-1) +6\chi ({\cal O}_S),
\end{equation}
\begin{equation}
\label{K3}
K_X^3=-\frac{9}{2}d^2+ \frac{27}{2}d+gd +18(g-1)-
30\chi ({\cal O}_S)-24\chi ({\cal O}_X),
\end{equation}
\begin{equation}
\label{c2L}
c_2\cdot L= -\frac{1}{4}d^2 +\frac{5}{2}d +2(g-1)+6 \chi ({\cal O}_S)
\end{equation}
\begin{equation}
\label{c3}
c_3=\frac{1}{4}d^2- \frac{1}{2}d -10(g-1)+gd +
24\chi ({\cal O}_S)-30\chi ({\cal O}_X).
\end{equation}
To prove (\ref{KL2}) we use the genus formula. (\ref{K2L}) follows from
\ci{a-s}, Proposition 2.1, after having realized that
$K_X^2\cdot L=$ ${K_X}^2_{|S}=$  $(K_S-L_{|S})^2$. The formula
for $K_X^3$ follows by ``cutting" (\ref{deg2dpf5}) with
$K_X$ and by using the above expressions for
$K_X\cdot L^2$, $K_X^2\cdot L$, and the fact that, by Hirzebruch-
Riemann-Roch on a threefold, $c_1c_2=24 \chi ({\cal O}_X)$.
The proof of (\ref{c2L}) is similar. (\ref{c3}) is obtained from
(\ref{deg3dpf5}) by first plugging the expression for
$(c_1^2-c_2)$ that one gets form (\ref{deg2dpf5})
and then by plugging the above relations into it. 

\smallskip
Finally we record the expression for the Hilbert polynomial of a 
threefold 
$X \subseteq \Q{5}$:
\begin{equation}
\label{chioxt}
\chi  ({\cal O}_X (t))  =   \frac{1}{6}d\,t^3+[\frac{1}{2}-
\frac{1}{2}(g-1)]\,t^2+     
[\frac{1}{3}d - \frac{1}{2}(g-1) +\chi ({\cal O}_S)]\,t +\chi ({\cal 
O}_X),   
\end{equation} 
which is an easy consequence of Hirzebruch-Riemann-Roch on a 
threefold  (cf. \ci{ha}, page 437) and  of  the formul\ae \ 
above.

\begin{fact}
\label{unirat}
{\rm (Unirationality of the Hilbert scheme)
 Let $H$
 the connected component of the Hilbert scheme of $\Q{n}$ containing 
the point
corresponding to  a fixed $X\subseteq \Q{5}$. Denote by $\frak H$ 
the open subscheme of $H$
 corresponding to nonsingular subvarieties. Assume that
every subvariety, $X'\in {\frak H}$, admits a resolution of its 
ideal sheaf
of the following form:
$$
0 \to \odix{\Q{n}}^s \to E \to {\cal I}_{X'}(c_1(E)) \to 0,
$$
where $E$ is a fixed,  locally free  sheaf independent of $X'$
and $s$ is a fixed positive integer.

\noindent
Under the above assumptions $\frak H_{red}$ is integral and unirational.
In fact it is enough to observe that
the natural  rational map 
${\Bbb P}(\wedge^s H^0(E)^{\vee}) --> {\frak H}$ is a dominant one.
}
\end{fact}

\medskip
 In the present  context,  Lemma 2.3 of \ci{a-s} gives the following:

\begin{fact}
\label{dimensionhilbert}
{\rm 
(Smoothness and dimension of the Hilbert Scheme)
}
{\rm If $h^i(\odix{X})=0$, $i\geq 1$, then $h^i({\cal N}_{X,\Q{n}})=0$,
$\forall i\geq 1$,
$\frak H$ is nonsingular  and of
dimension $h^0({\cal N}_{X,\Q{n}})$ at $X$. Riemann-Roch on a threefold,
$n_1=K_X+5L$, $n_2=(d/2)L^2$, and formul\ae \ 
(\ref{KL2}), (\ref{K2L}), (\ref{c2L}) give, for $n=5$:
$$
\chi ( {\cal N}_{X,\Q{n}})= -\frac{5}{4}d^2 + 10d + 10(g-1) + 
5\chi (\odix{S}).
$$}
\end{fact}

\begin{fact}
\label{hilbertforci}
{\rm  
(The Hilbert scheme of complete intersections) If $X\subseteq \Q{n}$ 
is a complete intersection of type $(2,i,j)$ in $\pn{n+1}$, then
the corresponding Hilbert scheme $\frak H$ is  integral, nonsingular 
nd rational. 

\noindent
For $i<j$:
\begin{eqnarray}
\dim {\frak H}=P(n;i,j)\hspace{-0.1in}&:=&  [  
B(n+1+i,n+1)-B(n+1+i-2,n+1)-1 ]
+  \nonumber  \\
&\ & [  B(n+j,n) - B(n+j-2,n) -1],   \nonumber
\end{eqnarray}
where $B(a,b):=a!/[b!(a-b)!]$, is the usual binomial coefficient.

\noindent
For $i=j$: 
$$
\dim {\frak H}=Q(n;i):=2[ B(n+1+i,n+1)-B(n+1+i-2,n+1) -2 ].
$$
}
\end{fact}

\bigskip
The following gives:  1) a method to construct codimension 
two subvarieties of $\Q{n}$ using vector bundles;  2) a way to reconstruct
the ideal sheaf of a codimension two subvariety, $X$, given enough sections
of twists of its dualizing sheaf $K_X$. 

\begin{fact}
\label{hartshorne}
{\rm The following is a Bertini-type Theorem  due to Kleiman; see \ci{kl}.
 Let $ E,$ $ F$ two  vector bundles on $\Q{n}$ of rank
$m$ and $m'$ respectively, such that
$E^{\vee}\otimes F$ is generated by its global sections. Let $\phi: E \to 
F$ 
be
an element of $H^0( E^{\vee}\otimes F)$. Define
$D_k(\phi)$ to be the closed subscheme of $\Q{n}$ defined, locally,
by the vanishing of
the $(k+1)\times (k+1)$ minors of a matrix representing  $\phi$. For
the general $\phi$ and for every $k$:

\noindent
 {\em 
{\rm a)} either $D_k$ is empty or it has
codimension $(m-k)(m'-k)$ and 
$D_k(\phi)_{sing}\subseteq D_{k-1}(\phi)$; in particular,
for $ n< (m-k+1)(m'-k+1)$, $D_k(\phi)$ is nonsingular;

\noindent
{\rm b)} for  $n\geq 5$, assuming that $D_k(\phi)$ has codimension two,
$D_k(\phi)$ is connected {\rm (see the remarks following Theorem 2.2 of
\ci{so-vdv})}. 
}

\noindent
The following fact, proved by Vogelaar, stems from an idea of Serre's;  
see  \ci{o-s-s}, Theorem I.6.4.2
or \ci{a-s} \S 2.3.
Let $X\subseteq \Q{n}$ be a local complete intersection of codimension 
two and
$a$ an integer such that the twist
$\omega_X(a)$ is generated by $s$ of its global sections. Then we have 
an exact sequence
$$
0 \to \odix{\Q{n}}^s  \to  F \to {\cal I}_X(n-a) \to  0,
$$
with $F$ locally free.
} 
\end{fact}

\subsection{A lifting criterion
and bounds for the genera of curves on $\Q{3}$}
\label{rothlifting}

The following is well known when $\Q{n}$ is replaced by
$\pn{n}$, see  \ci{boss} for example. The case
of $\Q{4}$ is proved in \ci{a-s}, Lemma 6.1. The general case
can be proved in the same way.  We used it as a tool to prove the 
finiteness of 
the number of families of nonsingular
codimension two subvarieties of $\Q{5}$ not of general type. See
\ci{bounded}, where we prove a more general statement.

\begin{pr}
\label{roth}
Let $X$ be an integral subscheme of degree $d$ and codimension two on 
$\Q{n}$,
 $n\geq 4$. Assume that for the general hyperplane section
$Y$ of $X$ we have
$
h^0({\cal I}_{Y,\Q{n-1}}(\s))\not= 0,
$
for some positive integer $\s$ such that $d>2{\s}^2$.
Then\,
$
h^0({\cal I}_{X,\Q{n}}({\s}))\not= 0.
$
\end{pr}

\medskip
\begin{evidence}
\label{1.1}
{\rm $C$ is an integral curve lying on a smooth three-dimensional quadric
$\Q{3}$, 
$k$ is  a positive integer,
$S_k$ is   an integral surface in
$|\odixl{\Q{3}}{k}|$
containing
$C$, 
$d$ 
and
$g$
 are the degree and the geometric genus of
$C$, respectively.}
\end{evidence}

\begin{defin}
\label{1.2}
{\rm Define
$n_0$ and   $\epsilon$
when 
$d>2k(k-1)$
and
$\ \theta_0$ 
and 
$\epsilon'$
when
$d\leq 2k(k-1)$ 
as follows:}

$$
\begin{array}{lrllllcc}
 &n_0 :=&  \lfloor \frac{d-1}{2k}\rfloor+1, \qquad
&d \equiv&   -\epsilon \quad (mod \ 2k), 
\quad  &0 \leq \epsilon\leq 2k-1;& \\
& & & & & & \\
&{\theta}_0 :=& \lfloor \frac{d-1}{2k}\rfloor +1, \qquad
&d \equiv&  -\epsilon'\quad (mod\ 2\theta_0),\quad  &0  \leq
\epsilon'\leq 2\theta_0-1.&
\end{array}
$$
\end{defin}

The following class of curves plays a central role in the understanding 
of the curves whose genus is the maximum possible. Arithmetically
Cohen-Macaulay is denoted by a.C.M..

\begin{defin}
\label{1.3}
{\rm  A curve $C$ as in (\ref{1.1}) is said to be in the class
${\frak S} (d,k)$,
if it is nonsingular, projectively normal and linked, in a complete 
intersection
on 
$\Q{3}$ 
of type 
$(k, n_0)$ 
if 
$d>2k(k-1)$  
(of type $(\theta_0, k)$ 
if 
$d\leq 2k(k-1)$), 
to an  {\it a fortiori} a.C.M. curve 
$D_{\epsilon}$
($D_{\epsilon'}$, respectively)
of  degree
$\epsilon$
($\epsilon'$ 
respectively)
lying on a quadric surface 
hyperplane section of 
$\Q{3}.$}
\end{defin}

\begin{pr}
\label{1.4}
{\em (Cf. \ci{de}.)}
  Notation as in {\rm (\ref{1.1})}  and {\rm Definition \ref{1.2}}. 
Assume  first that
$d>2k(k-1)$. Then

\noindent
{\rm (a)}  

$$
g-1\leq \pi(d,k) -\Xi
$$
where

\[ \pi(d,k)=
\left\{
\begin{array}{ll}
\frac{d^2}{4k}
+\frac{1}{2}(k-3)d-\frac{\epsilon^2}{4k}-\epsilon(\frac{k-\epsilon}{2}),
 & \mbox{if  $\  0\leq \epsilon\leq k,$}  \\
\ \\
\frac{d^2}{4k}+\frac{1}{2}(k-3)d-(k-\tilde \epsilon)
(\frac{\tilde \epsilon}{2}-\frac{\tilde \epsilon}{4k} + \frac{1}{4}),
&  \mbox{if $\  k+1\leq  \epsilon \leq 2k-1,\  
\tilde \epsilon:=\epsilon-k;$} 
\end{array}
\right. \]
and

\[ \Xi=\Xi(d,k)=
\left\{
\begin{array}{ll}
0 & \mbox{ if $\ \epsilon=0,\ 1,\ 2,\  2k-1,$}  \\
\ \\
1 &  \mbox{ if  otherwise.}
\end{array}
\right. \]
{\rm (b)} The bound is sharp for
$\epsilon =0$, $1,$ $2,$ $3,$ $2k-2,$ $2k-1$.
A curve achieves such a maximum possible genus if and only if
it is in the class
${\frak S} (d,k)$, except, possibly, the cases 
$\epsilon=3$,  $2k-2$.

\medskip
\noindent
 Assume $d\leq 2k(k-1)$. Then statements {\rm a)} and {\rm b)}, with
$\pi'(d,k)=\pi(d, \lfloor \frac{d-1}{2k}\rfloor +1)= \pi(d, \theta_0)$ 
and with
 $\Xi',$
 $\epsilon'$,  
$(\theta_0,k)$
and
$D_{\epsilon'}$
replacing
$\Xi,$ 
$\epsilon,$
$(k,n_0)$
and
$D_{\epsilon}$,
 respectively, hold.
\end{pr}

\begin{cor}
\label{coreasybound}
{\em (See \ci{a-s}, Proposition 6.4 for the case $d> 2k(k-1)$.)}
Notation as above. Then
$$
g-1\leq \frac{d^2}{4k} + \frac{1}{2} (k-3)d.
$$
\end{cor}

\begin{pr}
\label{boundasep}
{\em (Cf. \ci{a-s}, Proposition $6.4$.)}
Let $C$ be an integral curve in $\Q{3}$, not contained in any
surface of degree strictly less than $2k$. Then:
$$
g-1\leq \frac{d^2}{2k} +\frac{1}{2}(k-4)d.
$$
\end{pr}

\subsection{An  inequality}
\label{useful<} 
In this section we prove an inequality which is an essential
tool for  our proof of the classification of scrolls over
surfaces on $\Q{5}$.

Let $X\subseteq \Q{5}$ be a three dimensional, nonsingular variety,
 $\s$ be the smallest integer for which there exists
 a hypersurface $V$ in  $|{\cal I}_{X,\Q{5}}(\s)|$ and
 $\cal N$ the normal bundle of $X$. By the minimality of $\s$, 
the natural section $\odix{X} \to \check{\cal N} (\s)$ is not the 
trivial one.
The transposed of this section defines the  sheaf of ideals of
$\odix{X}$ of the singular locus of $V$ restricted to $X$.
Let us denote by $\tilde \Sigma$ the associated scheme. We obtain 
the surjection
${\cal N} \to {\cal I}_{\tilde \Sigma}(\sigma)$.

\begin{defin}
\label{Divisor} 
{\rm Let $D$ be the  divisorial component of $\tilde \Sigma$, i.e.
the unique effective Cartier  divisor of $X$ whose sheaf of ideals is 
the smallest
sheaf of principal ideals containing ${\cal I}_{\tilde \Sigma}$; $D$
may be empty.
Let $\Sigma$ be the one  dimensional component of $\tilde \Sigma$, i.e.
the  scheme associated with the sheaf of ideals
${\cal I}_{\tilde \Sigma}(D)$; $\Sigma$ is either empty or of pure 
dimension one.}
\end{defin}

From the above we get that the following two facts hold.

\begin{fact}
\label{homdim1}
{\rm  the sheaf ${\cal I}_{\Sigma}$ is either $\odix{X}$ or it has 
homological dimension one.}
\end{fact}
 
\begin{fact}
\label{ISIGMAspanned}
{\rm ${\cal I}_{\Sigma}(\s L - D) = {\cal I}_{\tilde \Sigma}(\s L)$;   
in particular,
${\cal I}_{\Sigma}(\s L - D)$ is generated by global sections since it is 
a quotient of 
$\cal N$ 
which is a quotient sheaf  of the globally generated sheaf 
${\cal T}_{\Q{5}}$.
}
\end{fact}

\begin{pr}
\label{bunchof<>}
Let $s_i$, $i=1,$ $2,$ $3$ be the Segre classes of
${\cal I}_{\Sigma}(\s L -D)$. 
Then\, 
$
s_1s_2\geq s_3 \geq 0,
$
$s_1$ and $s_1^2-s_2$ are represented by effective cycles. 
Moreover,
\begin{equation}
\label{55}
\chi (\odix{S})\geq \frac{1}{6\s}[(d-12\s )(g-1) + (\frac{1}{4}\s + 
\frac{3}{2})d^2
-\frac{13}{2}\s d] -\frac{1}{6\s}[\frac{1}{2}dL^2 - (K_X+5L)^2]D.
\end{equation}

\end{pr}

\noindent
{\em Proof.}
 By  Fact \ref{ISIGMAspanned} there is a surjection
$\odix{X}^{m} \to $${\cal I}_{\Sigma}(\s L -D)$, for
some $m$. By Fact \ref{homdim1} the kernel, $F$, of this surjection
is locally free. By the definition of Segre classes, $s_i=c_i(\check F)$.
The first part of the proposition follows from \ci{boss}, Lemma 5.1.

\noindent
As to the proof of the last inequality,
first we compute the Chern classes  $C_i$ of 
${\cal I}_{\Sigma}(\s L -D)$ using the following exact sequence
which is the Koszul resolution of ${\cal I}_{\Sigma}(\s L -D)$:
$$
0 \to \odixl{X}{K_X -\s L + D} 
\to {\cal N} \to {\cal I}_{\Sigma}(\s L -D) \to 0;
$$
we get

\noindent
$C_1=\s L -D,$

\noindent
$C_2= \frac{1}{2} d L^2 - (K_X+5L)(\s L -D) + (\s L -D)^2,$

\noindent
$C_3=- \frac{1}{2}d L^2(K_X +5L)+  \frac{1}{2}dL^2(\s L -D) +
 (K_X +5L)^2(\s L -D)-
2(K_X +5L )(\s L- D)^2 + (\s L-D)^3.$

\noindent
The  Segre classes of ${\cal I}_{\Sigma}(\s L -D)$
are $s_1=C_1$, $s_2=C_1^2- C_2$, $s_3= C_1^3- 2C_1C_2 + C_3$.
We  make explicit these Segre classes using the formul\ae \ for
the $C_i$. Then we use (\ref{K2L}) and (\ref{KL2}). We  now 
 use the part of the proposition that we have just proved: 
 ({\ref{55}) is $s_3\geq 0$.
\blacksquare

\subsection{Special vector bundles on quadrics}
\label{vbundlesonq}
\begin{fact}
\label{spinorbundles}
{\rm (Spinor Bundles) Here we collect some properties of
spinor bundles on quadrics. See   \ci{a-ott}. 

\noindent
Let ${\cal S}$ be the spinor bundle on an odd-dimensional quadric and
${\cal S}'$, ${\cal S}''$ be the two spinor bundles on an 
even dimensional quadric;
 if $n$ is the dimension of the quadric, 
the rank of these bundles is $2^{ \lfloor \frac{n-1}{2} \rfloor}$.

\noindent
For $n=2m+1$ (odd)  we have an exact sequence:
$$
0 \to {\cal S} \to  \odix{\Q{n}}^{ 2^{m+1} } \to {\cal S}(1) \to 0;
$$
for $n=2m$ (even) we have exact sequences:  
$$
0 \to S \to \odix{\Q{n}}^{ 2^{m} } \to S(1) \to 0,
$$
where $S$ denotes either ${\cal S}'$ or ${\cal S}''$.

\noindent
For $n=2m+1$ we have  ${\cal S}^{\vee} \simeq {\cal S}(1)$.
For $n=4m$ we have  ${{\cal S}'}^{\vee} \simeq {\cal S}'(1)$ and
${{\cal S}''}^{\vee} \simeq {\cal S}''(1)$;
for $n=4m+2$ we have  ${{\cal S}'}^{\vee} \simeq {\cal S}''(1)$ and
${{\cal S}''}^{\vee} \simeq {\cal S}'(1)$.
Let $i: \Q{2k-1} \hookrightarrow \Q{2k}$ be a nonsingular hyperplane 
section;
then $i^* {\cal S}' \simeq$ $i^* {\cal S}'' \simeq$ $\cal S$.
Let $j: \Q{2h} \hookrightarrow \Q{2h+1}$ be a nonsingular hyperplane 
section;
then $j^* {\cal S} \simeq {\cal S}' \oplus {\cal S}''$.

\noindent
An analogue of Horrocks splitting criterion holds on quadrics;
recall that spinor 
bundles carry no intermediate cohomology:

\noindent
{\em let  $E$ be a  vector bundle  on $\Q{n}$ then $h^i(E(t))=0$, 
$0<i<n$, 
$\forall t \in \zed$ if and only if
$E$  splits as the direct sum of line bundles and
twists of spinor bundles of $\Q{n}$.}

\noindent
The Chern polynomial of ${\cal S}(l)$ on $\Q{5}$ is:
\begin{eqnarray}
c({\cal S}(l)) &=& 1 + (4l-2)h + (6l^2-6l+2)h^2 +(4l^3-6l^2+4l-1)h^3+ 
\nonumber \\
&+& (l^4-2l^2+2l^2-l)h^4. \nonumber
\end{eqnarray}
The Chern polynomial of ${\cal S}'(l)$ on $\Q{6}$ is:
\begin{eqnarray}
c({\cal S}'(l)) &=& 1+ (4l-2 )h+ (6l^{2}-6l+2
 )h^{2}+[ (4\,l^{3}-6l^{2}+4l
 )h^{3} -2{\Lambda_1}] +     \nonumber \\ 
&+&[ ( l^{4} -2l^{3}   +  2l^{2})h
^{4} -2l\,{\Lambda_1}h]     \nonumber
\end{eqnarray}
Replacing $\Lambda_1$ by $\Lambda_2$ in the formula above, 
we get $c({\cal S}''(l))$.
}
\end{fact}

\begin{fact}
\label{cayley}
{\rm (Cayley bundles) See \ci{ottcayley}.
 On $\Q{5}$ there is a family of rank two stable vector bundles, called
Cayley bundles. Each Cayley bundle $\cal C$ has 
Chern classes $c_1=-1$, $c_2=1$ and  ${\cal C}(2)$ is generated by
 global sections.
Every stable $2$-bundle on $\Q{5}$ with Chern classes $c_1=-1$, $c_2=1$ is
 a Cayley bundle. Cayley bundles are parameterized by 
 a fine moduli space isomorphic to
$\pn{7} \setminus \Q{6}$.
A Cayley bundle restricts, on a  $\Q{4}$, to a bundle of type 
$\dot{E}$  which appears in the description of
Type 10) of \ci{a-s}, page 44.
The Chern polynomial of a ${\cal C}(l)$ is:
$c({\cal C}(l))= 1 + (2l-1)h + (l^2-l+1)h^2.$ 
}

\end{fact}

\section{CLASSIFICATION FOR $d\leq 10$}
\label{d<=10}

\subsection{The list}
\label{list}
In what follows:

 - $((a,b,c),{ \cal O}(1))$ denotes the polarized pair given by a  
complete intersection
of type  $(a,b,c)$ in $\pn{n+1}$ and the restriction of the hyperplane 
bundle to it;

- $(X,L)$  denotes the polarized pair given by
a variety $X\subseteq \Q{n}$ and $L:= \odixl{\Q{n}}{1}_{|X}$; if we do not
explicitly say the contrary, the embeddings are projectively  normal
in $\pn{n+1}$; this fact  follows from the cohomology of the presentation
of the ideal sheaf;

- by a {\em presentation} of the ideal sheaf ${\cal I}_X$ we mean
an injection of  locally free coherent sheaves on $\Q{n}$, $\phi: 
 E \to F$, such that
$coker(\phi)\simeq {\cal I}_X(i)$, where $i=c_1(F)-c_1(E)$; we write 
the presentations
so that the integer $i$ is the smallest for which the sheaf $F$ is 
generated by global sections, so that
for that $i$ so will be the sheaf ${\cal I}_{X}(i)$;

 - $\frak H$ denotes the Hilbert scheme of $\Q{n}$ 
of a variety fixed by the context; see Fact \ref{unirat};

- $P(n;i,j)$ and $Q(n;i)$ are defined in Fact \ref{hilbertforci};

-  a digit ``$\#$)"  refers to the type of the surface section as in
\ci{a-s}, page 44 (where $d\leq 8$); type $Z^F_{10}$ refers to the paper
\ci{gross}

- $g,$ $q$ and  $p_g$  denote the sectional genus of the embedding
line bundle, the irregularity and geometric genus of a surface section, 
respectively.

\begin{tm}
\label{classificationd<12}
Let $X\subseteq \Q{n}$, $n\geq 5$, a codimension 
two nonsingular subvariety of degree $d\leq 10$. Then the pair
$(X,L)$, a presentation of the ideal of $X$ on $\Q{n}$ and the 
Hilbert scheme,
 $\frak H$,
of $X$ on $\Q{n}$ are as follows.
 
\medskip

\noindent
{\rm ($\bullet$)} $d=2$

\smallskip
\noindent
\underline{\rm Type A):}\quad  $((1,1,2), {\cal O}(1))$;\quad
$\odixl{\Q{n}}{-1} \to \odix{\Q{5}}^2$;\quad $\frak H$ is integral, 
nonsingular,
rational, of dimension $Q(n;1)$;\quad  {\rm  2)};\quad $g=q=p_g=0$.    

\medskip

\noindent
{\rm ($\bullet$)} $d=4$

\smallskip

\noindent
\underline{\rm Type B):}\quad  $((1,2,2),{\cal O}(1))$;\quad
$\odixl{\Q{n}}{-1} \to \odixl{\Q{n}}{1} \oplus \odix{\Q{n}}$;\quad 
$\frak H$ is integral, nonsingular, rational and of dimension
$P(n;1,2)$;\quad {\rm 6)};\quad $g=1$, $q=p_g=0$.

\smallskip
\noindent
\underline{\rm Type C):}\quad   $n=6$, $(\pn{1}\times \pn{3}, 
{\cal  O}(1,1))$;\quad
 $\odix{\Q{6}}^3
\to {\frak S}(1), $ with  ${\frak S}\simeq {\cal S}', {\cal S}''$;
\quad $\frak H$
consists of two connected components, which are both nonsingular, 
integral,
unirational and of dimension $15$;
\quad {\rm 5)};\quad $g=q=p_g=0$.

\smallskip
\noindent
\underline{\rm Type D):}\quad    $n=5$, 
$({\Bbb P}( {\odixl{\pn{1}}{1}}^2\oplus
\odixl{\pn{1}}{2}), \xi)$;\quad   
$\odix{\Q{6}}^3 \to {\cal S}(1);$ \quad $\frak H$ is integral, nonsingular,
 unirational
and of dimension $15$; \quad {\rm 5)};\quad $g=q=p_g=0$.

\medskip
\noindent
{\rm ($\bullet$)} $d=6$

\smallskip

\noindent
\underline{\rm Type E):}\quad   $((1,2,3),{\cal O}(1))$;\quad
$\odixl{\Q{n}}{-1} \to \odixl{\Q{n}}{2}\oplus \odix{\Q{n}}$;\quad
${\frak H}$ is integral, nonsingular, rational and of dimension
$P(n;1,3)$;\quad {\rm 12)};
\quad $g=4$, $q=0$, $p_g=1$.

\smallskip
\noindent
\underline{\rm Type F):}\quad  $n=5$, $({\Bbb P}( {\cal 
T}_{\pn{2}}),\xi)$, 
embedded using a general
codimension one linear system 
${\frak l}\subseteq| \xi_{ {\cal T}_{\pn{2}} }|$;\quad
$\odix{\Q{5}} \to {\cal C}(2)$;\quad
$\frak H$ is integral, nonsingular, unirational and of dimension
$20$;\quad {\rm 10)};\quad $g=1$, $q=p_g=0$.

\smallskip
\noindent
\underline{\rm Type G):}\quad
  $n=5$, $f: X \to \pn{1} \times \pn{2}=:Y$
a double cover, branched
along a divisor of type $\odixl{Y}{2,2}$,
$L\simeq   p^* \odixl{ Y }{1,1}$;\quad
${\odixl{\Q{5}}{-1}}^2 \to \odix{\Q{5}}^3$; $\frak H$ is 
integral, nonsingular, unirational and of dimension $30$;\quad
{\rm 11)};\quad $g=2$, $q=p_g=0$.

\medskip

\noindent
{\rm ($\bullet$)} $d=8$
\smallskip

\noindent
\underline{\rm Type H):}\quad  $((1,2,4), {\cal O}(1))$;\quad
$\odixl{\Q{n}}{-1} \to \odixl{\Q{n}}{3} \oplus \odix{\Q{n}}$;\quad
$\frak H$ is integral, nonsingular, rational and of dimension
$P(n;1,4)$;\quad {\rm 20)};\quad $g=9$, $q=0$, $p_g=5$.

\smallskip
\noindent
\underline{\rm Type I):}  $((2,2,2), {\cal O}(1) )$;\quad 
$\odixl{\Q{n}}{-2} \to \odix{\Q{n}}^2$;\quad
$\frak H$ is integral, nonsingular, rational and of dimension
$Q(n;2)$;\quad {\rm 19)};\quad $g=5$, $q=0$, $p_g=1$.

\smallskip
\noindent
\underline{\rm Type L):} 
$n=5$, $({\Bbb P}( E), \xi)$, $E$ a rank two vector bundle on $\Q{2}$ as in
{\rm \ci{io3}};\quad $\odix{\Q{5}}^4 \to {\cal S}(1)\oplus 
\odixl{\Q{5}}{1}$;\quad
$\frak H$ is integral, nonsingular, unirational and of dimension
$35$;\quad {\rm 18)};\quad $g=4$, $q=p_g=0$.

\medskip
\noindent
{\rm  ($\bullet$)} $d=10$
\smallskip

\smallskip
\noindent
\underline{\rm Type M):}  $((1,2,5), {\cal O}(1) )$;\quad
$\odixl{\Q{n}}{-1} \to \odixl{\Q{n}}{4} \oplus \odix{\Q{n}}$;\quad
$\frak H$ is integral, nonsingular, rational and of dimension
$P(n;1,5)$;\quad $g=16,$ $q=0,$ $p_g=14$.

\smallskip
\noindent
\underline{\rm Type N):}\quad $n=5$, $f_{|K_X+L|}: X \to \pn{1}$ 
is a fibration
with general fiber a Del Pezzo surface $F$, $K_F^2=4$, $K_X=-L+
f^*\odixl{\pn{1}}{1}$;\quad $\odixl{\Q{5}}{-1}^2 \to 
\odixl{\Q{5}}{1} \oplus \odix{\Q{5}}^2$;\quad
$\frak H$ is integral, nonsingular, unirational and of dimension 
 $60$; type $Z_{10}^F$;\quad $g=8$, $q=0$, $p_g=2$. 

 \end{tm}

\begin{rmk}
{\rm
In this remark, by the symbol ${\rm Q} \stackrel{(a,b)}{\sim} 
{\rm R}$, we mean
that {\em every} variety of Type Q) is linked to {\em a} variety of Type R)
in a complete intersection of type $(a,b)$ on $\Q{5}$. Using
Lemma \ref{pesk} and the presentations of the ideals of the varieties 
of the above theorem
we see that:
\noindent
${\rm A} \stackrel{(1,2)}{\sim} {\rm A}$, 
${\rm A} \stackrel{(1,3)}{\sim} {\rm B}$,
${\rm A} \stackrel{(1,4)}{\sim} {\rm E}$, 
${\rm A} \stackrel{(2,2)}{\sim} {\rm G}$,
${\rm A} \stackrel{(1,5)}{\sim} {\rm H}$, 
${\rm A} \stackrel{(1,6)}{\sim} {\rm M}$,
${\rm A} \stackrel{(2,3)}{\sim} {\rm N}$, 
${\rm B} \stackrel{(2,2)}{\sim} {\rm B}$, 
 ${\rm B} \stackrel{(2,3)}{\sim} {\rm I}$, 
${\rm G} \stackrel{(2,2)}{\sim} {\rm A}$,
${\rm G} \stackrel{(2,3)}{\sim} {\rm G}$, 
${\rm G} \stackrel{(2,4)}{\sim} {\rm N}$,
${\rm I} \stackrel{(2,3)}{\sim} {\rm B}$, 
${\rm I} \stackrel{(2,4)}{\sim} {\rm I}$,
${\rm I} \stackrel{(3,3)}{\sim} {\rm N}$, 
${\rm N} \stackrel{(3,3)}{\sim} {\rm I}$.

\noindent
The simple details are left to the reader. As for Type F), see Proposition
\ref{example}.

}
\end{rmk}

\subsection{The proof}
\noindent
{\em Proof of Theorem} \ref{classificationd<12}. The degree $d$ 
is always an even integer
by
Remark \ref{rmkdiseven}. The statements of the Theorem concerning
complete intersections follow from
Fact \ref{hilbertforci} and \ci{ha}, III.9.Ex. 9.6. 
 In the sequel, we do not deal with complete intersections.

\medskip

\noindent
{\em Claim. The only nonsingular surfaces on $\Q{4}$ which can 
be a general hyperplane
section of a threefold on $\Q{5}$ of degree $d\leq 10$ are: types}
5), 10), 11) {\em and} 18) {\em from} \ci{a-s} {\em and type} $Z^{10}_F$ 
{\em from} \ci{gross}.
}

\medskip
\noindent
{\em Proof of the Claim.}
Let $d\leq 8$.
\ci{a-s} page 44 contains the complete list of nonsingular surfaces on 
a $\Q{4}$
of degree $d\leq 8$.  Not all of them can be a general hyperplane
section of a threefold on  $\Q{5}$. The complete list
of linearly normal, nonsingular subvarieties of projective space of degree
$d\leq 8$ is given in \ci{io1} and \ci{io3}. We are going to use
these results jointly. 

 \noindent
Since $H^2(\Q{5},\zed )\simeq \zed <h^2>$, the surface section, $S$, 
of a degree $d$ threefold $X\subseteq \Q{5}$ has cohomology class
$[S]=(d/2) \Lambda_1 + (d/2) \Lambda_2$.
This implies that the surfaces of type 
1), 3), 4), 7), 8), 13), 14), 15) and 16) of \ci{a-s} page 44 cannot be 
nonsingular hyperplane sections of any threefold on $\Q{5}$.
\noindent
Types 2), 6), 12), 19) and 20)  are complete intersections and we drop
 them from the list.

\noindent
In what follows, assume that $S$ is a surface of a given type and that
$X\subseteq \Q{5}$ is a threefold with general surface section $S$. 
We now exclude
types 9) and 17).
Type 9). If  $X$ existed,  
a comparison with Ionescu's list
\ci{io1}  would force $X$ to be a rational scroll
over a curve contradicting Proposition
\ref{pbundlesovercurves}.
Type 17). This type has 
 sectional genus $g=3$ so that,
 according to Ionescu's list \ci{io3}, 
 $X$ would have to be either  a scroll over an
elliptic curve, 
a scroll over $\pn{2}$ or $X$ would have to admit a morphism onto 
$\pn{1}$ with
all fibers quadric surfaces. We exclude 
the first case  because type 17) is simply connected, the second one
 by Proposition \ref{baseofscroll}. The last one
 would imply, after having cut (\ref{deg2dpf5}) with a general fiber
$F\simeq \Q{2}$,  the contradiction $d=6$.

\noindent
It follows that, except for the case of complete intersections,
 only the following types are admissible as surface
sections of codimension two nonsingular subvarieties
of quadrics when $d\leq 8$: \,5), 10), 11), 18).

\noindent
Let $d=10$. We  employ the same technique
as above using \ci{gross} and \ci{fa-li} instead of \ci{a-s} and 
\ci{io1}, \ci{io3}.
Looking at the list in \ci{gross} we exclude  cases $Z^{10}_A$ and 
$Z^{10}_B$
since they do not have a 
balanced cohomology  class. Cases $Z^{10}_D$ and $Z^{10}_E$
cannot occur  by  \ci{fa-li},  since they have sectional genus $g=7$.

\noindent
The case $C^{10}_A$, where  the sectional genus $g=4$ and 
the irregularity $q=1$ is excluded since, by
 \ci{fa-li}, we would have $q=0$. The cases
$C^{10}_B$ and $C^{10}_C$  are excluded in a similar way.

\noindent
The case $Z^{10}_C$, which is a rational surface
with $g=6$ is excluded as follows. According to \ci{fa-li}
there are only two types of threefolds of degree $d=10$ with sectional 
genus $g=6$;
the first is a Mukai manifold (i.e. $K_X=-L$), the second one a scroll
 over $\pn{2}$.
In the former case the surface section would have trivial canonical
 bundle, contradicting
its being  rational. 
The latter case is excluded by Proposition \ref{baseofscroll}. 

\noindent
The proof of the  Claim  is complete.

\medskip
\noindent 
We now show that all the  types of the claim occur as nonsingular 
surface sections
of threefolds on $\Q{5}$, that type $5$) is the only one that
can occur as a section of a fourfold
on $\Q{6}$ and that none of these types can occur as a surface 
section of any
$(n-2)$-fold on $\Q{n}$,
for $n\geq 7$.

\medskip
\noindent
{\em The case of Type $5)$.}

\noindent
Assume that $X\subseteq \Q{6}$ is a fourfold with surface section of type
$5$).  By Swinnerton-Dyer's classification of varieties of degree
$d=4$ (see \ci{io1} for example), $(X,L)$  is  of Type C); such a 
type occurs as a 
subvariety of $\Q{6}$
as pointed out in Proposition \ref{pbundlesovercurves}. 
$K_X=\odixl{X}{-2,-4}$, so that
$K_X(4)$ is generated by three global sections. 
Fact \ref{hartshorne} gives us the following exact sequence:
\begin{equation}
\label{F}
0 \to \odix{\Q{6}}^3 \to F \to {\cal I}_X(2) \to 0,
\end{equation}
where $F$ is locally free. We want to prove that
$F$ is isomorphic to either ${{\cal S}'}_6(1)$ or to ${{\cal S}''}_6(1)$,
where the subindeces refer to the fact that that the bundles are the 
spinor bundles of 
$\Q{6}$.
Consider a general threefold section $T$, which is of
type D), and a general surface  section $S\subset T$.
 We have $K_T \simeq K_X(1)\otimes \odix{T}$ and 
$K_S \simeq K_X(2)\otimes
\odix{S}$; there is a canonical identification between
$H^0(K_X (4))$, $H^0(K_T(3))$ and
 $H^0(K_S(2))$ so that the bundles $\cal F_T$ and $\cal F_S$,
 that we obtain
repeating for $T$ and $S$ the construction we have done for $X$ using 
Fact
\ref{hartshorne}, satisfy
${\cal F}_T\simeq F_{|T}$ and ${\cal F}_S\simeq F_{|S}$.

\noindent
We know, from \ci{a-s}, that, on $\Q{4}$, 
${\cal }F_{|S}\simeq 
{{\cal S}'(1)}_4 \oplus {{\cal S}''}_4(1)$. Recall
that spinor bundles have no intermediate cohomology.
Let us  look at the long cohomology sequences associated with the exact 
sequences:
$$
0 \to {\cal F}_T(-1+t) \to {\cal F}_T(t) \to {\cal F}_S(t) 
\to 0, \quad t\in \zed.
$$
Firstly we deduce that $H^1({\cal F}_T(-\tau+t))$ surjects
onto $ H^1({\cal F}_T(t))$, for every fixed $t$ and every $\tau \geq 0$;
 Serre Duality and Serre Vanishing imply that $h^1({\cal F}_T(t))=0$, 
$\forall t$.
The vanishing of $h^i({\cal F}_T(t))$ for $2\leq i \leq 4$ are 
dealt with similarly.
We have proved that ${\cal F}_T(t)$ has no intermediate cohomology, 
so that,
by the analogue of Horrocks criterion in Fact
\ref{spinorbundles}, ${\cal F}_T$ splits as a direct sum
of line bundles and twists of spinor bundles on $\Q{5}$.
Since the rank of $\cal S$ is four we see that 
either ${\cal F}_T\simeq {\cal S}(j)$ or it splits completely as the
direct sum of line bundles ${\cal F}_T\simeq 
\bigoplus_{i=1}^4 \odixl{\Q{5}}{a_i}$.
Using  the Castelnuovo-Mumford $0$-regularity criterion for  global

 generation (see
 \ci{a-s}, where it is proved for sheaves on $\Q{4}$; 
the case of any $\Q{n}$
is analogous)
we  see that ${\cal F}_T$ is generated by global sections as soon as
$h^5({\cal F}_T(-5))=0$ which follows from the cohomology sequence
 associated with the sequence (\ref{F}) twisted by $-5$ once we observe that
$h^4(\Il{T}{\Q{5}}{-3})=3$ and $h^5(\Il{T}{\Q{5}}{-3})=0$. 
Recall that
$1$ is the smallest integer $j$
for which the spinor bundles twisted by $j$ are generated by global sections.
It follows
that, for the splitting type of ${\cal F}_T$, we have either
$j\geq 1$   
or $a_i\geq 0$, $\forall i$. 

We can compute the Chern classes of ${\cal F}_T$ using (\ref{F}), 
the invariants
of \,$T$ and 
Lemma \ref{lmccis}. Comparing Chern polynomials we deduce that
${\cal F}_T$ cannot split as the direct sum of line bundles, and that, 
once it is a twist
of the spinor bundle ${\cal S}$, $j=1$: ${\cal F}_T \simeq {\cal S}(1)$.

\noindent
We repeat the argument, replacing $S$ with $T$ and $T$ with $X$, to see
that either $F\simeq {\cal S}'(1)$  or $F\simeq {\cal S}''(1)$.

\noindent
We  have proved that every fourfold
on
$\Q{6}$
 with surface section of type 5) is as in Type C) and has the prescribed
presentation for its ideal sheaf;
we have also proved that 
 {\em every} threefold on $\Q{5}$ with surface section
of type
5) is of Type  D)  and has the prescribed
presentation for its ideal sheaf.  Conversely, since ${\cal S}(1)$, 
${\cal S}'(1)$
and ${\cal S}''(1)$ are globally generated, we use Fact \ref{hartshorne}
and our maximal list of varieties of degree $d=4$ to prove that the variety
$D_2(\phi)$ is
 as in C) or D), where $\phi$ is a general element of
$H^0(S(1)^3)$ and $S$ one of the three spinor bundles in question.
To be precise, Fact \ref{hartshorne}
implies that,  for a general $\phi$ on $\Q{6}$, 
$D_1(\phi)$  is either empty or has the expected codimension  six
and  $D_2(\phi)$ will be nonsingular outside
$D_1(\phi)$. Porteous' formula, \ci{a-c-g-h}, II.4.2 gives 
$[D_1(\phi)]= c_3(S(1))^2 - c_2(S(1))c_4(S(1))=0$; the Chern classes of
the spinor bundles are  listed in \ref{spinorbundles}.
It follows that $D_1(\phi)=\emptyset$.

\noindent
Using facts (\ref{unirat}) and (\ref{dimensionhilbert}) we conclude 
the proof for
type D). To complete the proof for type C) we remark that
${\cal S}'(1)$ distinguishes, via the choice of three general sections,
a nonsingular, integral
 component, say ${\frak H}'$, of $\frak H$. The same is true for
${\cal S}''(1)$ which defines another, distinct,  nonsingular component
${\frak H}''$. Since $\frak H$ is nonsingular,
${\frak H}= {\frak H}' \bigsqcup {\frak H}''$.

\noindent
The dimension of the two components, which are abstractly isomorphic, 
can be computed
using Riemann-Roch and Fact \ref{dimensionhilbert}.
The fourfold C) cannot be the hyperplane section of a fivefold; see \ci{io1}.

\medskip
\noindent
{\em The case of Type} 10).

\noindent
Assume that $X$ is a threefold on $\Q{5}$ whose general surface section, 
$S$, is of type
10). Since $K_S=-L_{|S}$ and 
the natural map $Pic(X)\to Pic(S)$ is injective by
Lefschetz theorem on hyperplane sections, we have $K_X=-2L$. Looking at
\ci{io1} for degree $d=6$ we see that either $(X,L)$ is as in Type F)
or $X$ is a scroll over $\Q{2}$; the latter case is not possible by 
Proposition
\ref{baseofscroll}.
We have $K_X=-2L$. 
 Fact \ref{hartshorne} yields an extension:
$$
0 \to \odix{\Q{5}} \to F \to {\cal I}_X(3) \to 0,
$$
with $F$ locally free of rank two.
By \ci{o-s-s}, 2.1.5 one shows that $F(-2)$ is stable. Using
Lemma \ref{lmccis} we deduce that, for the Chern classes of $F$, 
$c_1(F(-2))=-1$
and 
$c_2(F(-2))=1$. By Fact \ref{cayley} we see that $F(-2)$
is a Cayley bundle. Conversely, \ci{ottcayley}, Theorem 3.7 
ensures that the general section
of a normalized Cayley bundle twisted by $\odixl{\Q{5}}{2}$ vanishes
exactly along a variety of type F). As in \ci{a-ss}, page 209, we see
that our scrolls are parameterized by an open dense set, ${\frak U}$
 of a projective bundle
over the fine moduli space, $\pn{7}\setminus \Q{6}$, of 
these Cayley bundles.
This space is clearly rational and it has dimension 20.

\noindent
${\frak U}$ admits a natural morphism onto the Hilbert scheme
$\frak H$ of our scrolls. This morphism is one to one. To conclude
it is enough to observe that $\frak H$ is nonsingular by Fact
\ref{dimensionhilbert}, for then
the morphism in question  is  an isomorphism
by Zariski Main Theorem. 

\noindent
Note that, again by \ci{io1}, this threefold is the hyperplane section 
of only one fourfold, $\pn{2}\times \pn{2}$ embedded via the Segre
 embedding;
 this latter  can be projected smoothly to $\pn{7}$ but, after this 
embedding,
it does not lie on a smooth quadric $\Q{6}$ by Proposition \ref{nop2inq6}.

\medskip

\noindent
{\em The cases} 11), 18) {\em and} $ Z^{10}_F$.

\noindent
These cases are analogous to the one of type 5).
If the general surface section, $S$,  of a threefold $X$ on $\Q{5}$ 
is of type 11) then there is a morphism
with connected fibers $f:S\to \pn{1}$ all fibers of which are conics
and $ K_S(1)= f^*\odixl{\pn{1}}{1}$, so that the former sheaf 
 is generated by two global sections.

\smallskip
\noindent
{\em Claim.
$K_X(2)$ is generated by its global sections.}
Fix any point $x \in X$. Take any nonsingular hyperplane section
$S$ of $X$ through $x$; there are plenty of them since the dual variety 
$\hat X$
does not contain hyperplanes.
 Kodaira Vanishing implies that $H^0(X, K_X(2))$ surjects onto
$H^0(S, K_S(1))$ which in turn generates 
the stalk of $K_{S}(1)$. The claim follows.

\smallskip
\noindent
A computation analogous to the one of type 5) allows us to conclude that
if the threefold in question exists than it is of type G). To prove 
its existence,
we use Fact \ref{hartshorne} for a general morphism
$\phi: \odix{\Q{5}}^2 \to \odixl{\Q{5}}{1}^3$.
 This threefold cannot be
the hyperplane section of any fourfold by \ci{io1}.

\smallskip
\noindent
If the type of the general surface section is 18) then we have
a morphism
 $f:S \to \Q{2}$ which is the
blowing up of $\Q{2}$ at 10 points and  $K_S(2)\simeq f^*\odixl{\Q{2}}{1}$
is generated by four global sections. We argue as above and get
Type L).

\smallskip
\noindent
If X is a threefold on $\Q{5}$ with general section $S$ of type
$Z^{10}_F$ then the sectional genus $g=8$. By looking at the list
in \ci{fa-li}, we see that $X$ is a Del Pezzo fibration
$f:X \to Y$ over a curve $Y$ and that
$X$ is not the hyperplane section of any fourfold.  By looking at the proof
of \ci{fa-li} Proposition 4.1 we see that the base of the fibration, $Y$, 
is a rational curve
and  that $K_X(1) \simeq f^* \odixl{\pn{1}}{1}$ is generated by two 
global 
 sections. We argue as above and conclude that the type is N).
\blacksquare

\section{SCROLLS ON QUADRICS}
\label{secpbundles}

\subsection{Statement of the main result}
\label{main}

In this section we classify scrolls 
 as codimension two subvarieties of $\Q{n}$, for $n\geq 5$.
A {\em scroll}, here, is a nonsingular subvariety 
$X\subseteq  \Q{n} \subseteq \pn{n+1}$ which
admits a surjective morphism
$p:X\to Y$ to a lower dimensional variety $Y,$ such that $p$ has 
equidimensional fibers
 and the general scheme theoretic fiber is a linear subspace of 
$\pn{n+1}$ 
of the appropriate
dimension. 
 The case $\dim Y=0$ is the theory of maximal dimensional linear spaces 
in quadrics, 
a well known subject; see \ci{h-p}. 
From now on we assume $\dim Y>0$. 

\noindent
By standard arguments,
see \ci{fujitasemipositive} 2.7, we can assume, without loss of 
generality, that
$Y$ is nonsingular and that the polarized pair $(X,L) \simeq   
(\Bbb P ({\cal E}), \xi_{\cal E})$,
 where $L$ is the restriction of
the hyperplane bundle to $X$ and  ${\cal E}:= p_*L$ is
 an ample, rank $\mu:=\dim X -\dim Y +1$, locally free sheaf 
generated by 
its global sections.

\begin{rmk}
\label{gold}
{\rm
We  assume  that $n\geq 5$ since
surfaces on $\Q{4}$ which are  scrolls over curves have been 
classified
by Goldstein  in \ci{go}. They correspond to the surfaces of type
2), 3), 5), and 9)  in \ci{a-s}.
} 
\end{rmk}

In what follows the Types C), D), F) and I) below refer to Theorem
\ref{classificationd<12}; we say that a nonsingular threefold, $X$,
on $\Q{5}$ is of Type O), if it  has degree $d=12$ and it
is a scroll over a minimal $K3$ surface.

\begin{tm}
\label{maintm}
The following is the complete list of \,nonsingular codimension 
two subvarieties
of quadrics $\Q{n}$, $n\geq 5$, which are scrolls.

\noindent
{\rm Type C)}, $n=6$, $d=4$,  scroll over $\pn{1}$ and over $\pn{3}$;
 
\noindent
{\rm Type D)}, $n=5$, $d=4$, scroll over $\pn{1}$; 

\noindent
{\rm Type F)}, $n=5$, $d=6$, scroll over $\pn{2}$;

 \noindent
{\rm Type \,L)}, $n=5$, $d=8$, scroll over $\Q{2}$;

\noindent
{\rm Type O)},  $n=5$, $d=12$, scroll over a minimal $K3$ surface.
\end{tm}

\noindent
{\em Proof.}
The proof is the consequence of the  lengthy analysis that 
constitutes the rest of the paper. Here we give the reader  directions 
toward
the various relevant statements. 

\noindent
By    Fact \ref{n<7} we need to deal only with
the cases $n=5,\,6$.

\noindent
Scrolls over curves are classified by
Proposition \ref{pbundlesovercurves}.   They correspond
to types C) and D).

\noindent
The are no fourfolds  which are scrolls over a surface, by   
Proposition \ref{nop2inq6}. 

\noindent
The only fourfold which is a scroll over a threefold  is of type C), by 
Proposition \ref{p1bundlesinq6}.

\noindent
Threefolds which are scrolls over a surface have degree $d=6,\, 8,\, 12$
 by Proposition \ref{d=6812}  and the base surfaces are as in Proposition
\ref{baseofscroll}. The classification in degrees
$d=6,\, 8$ is complete; correspondingly we get
types F) and L).

\noindent
For an example and for the general properties of varieties  of Type O)
see Section \ref{d=12}).
\blacksquare

\subsection{Preliminary facts}
\label{prelfacts}
\medskip
 The Barth-Larsen theorem   implies that 
if $X$ is a nonsingular codimension two subvariety of 
$\Q{n}$, then the fundamental group $\pi_1(X)$ is trivial for
$n\geq 6$, and $Pic(X)\simeq \zed$, generated by the hyperplane bundle,
for
$n\geq 7$; see \ci{ba}.
Since  $Pic(X)\simeq \zed$ as soon as
$n\geq 7$,
 we have:

\begin{fact}
\label{n<7}
{\rm 
There are no codimension two scrolls on
 $\Q{n}$ for  $n\geq 7$ and, for $n=6$, any such is simply
connected.}
\end{fact}

It is therefore enough to study threefolds on $\Q{5}$ which are scrolls
over curves and  surfaces and fourfolds  on $\Q{6}$ which are
 scrolls over curves, surfaces, and threefolds. 

\smallskip
Let us begin the analysis by fixing some notation. We start with a
scroll  of degree $d$;
 let
$e_i:=c_i({\cal E})$, $x_i:=c_i(X)$ and  $y_i:=c_i(Y)$.
Since $p^*$ is injective it is harmless to denote
$p^* \alpha$ simply by $\alpha$ while performing computations in 
the cohomology ring of $X$.

\noindent
The tautological relation is 
\begin{equation}
\label{tautological}
\sum_{i=0}^{\mu} (-1)^i L^{\mu -i}\cdot e_i=0.
\end{equation}

\noindent
Finally, recall the usual exact sequence:
\begin{equation}
\label{euler}
0\rightarrow \odix{X}\rightarrow p^*({\cal E}^{\lor})\otimes L
 \rightarrow {\cal T}_X \rightarrow
p^*{\cal T}_{Y} \rightarrow 0,
\end{equation}
which is obtained by pasting together the relative Euler sequence
 \ci{ha}, II.8.13 and the short exact sequence 
associated with the epimorphism
$dp:{\cal T}_X\to p^*{\cal T}_Y$.

\subsection{Scrolls over curves on $\Q{5}$ and on $\Q{6}$}
\label{overcurves}

The following is proved independently of Theorem \ref{classificationd<12}.

\begin{pr}
\label{pbundlesovercurves}
Let $(X,L)$ be scroll over a nonsingular
curve $Y$, on $\Q{n}$.
 Then $(X,L)$ is one of the following:

\noindent
{\rm (\ref{pbundlesovercurves}.1)}
$ n=6,$  
 $ ({\Bbb P}_{\pn{1}}(\cal E),\xi_{\cal E})$, 
${\cal E}:={{\cal O}_{{\Bbb P}^1}(1)}^4.$

\noindent
{\rm (\ref{pbundlesovercurves}.2)}
$n=5,$ 
$({\Bbb P}_{\pn{1}}(\cal E),\xi_{\cal E})$, 
${\cal E}:={\cal O}_{{\Bbb P}^1}(2)
\oplus  {\odixl{\pn{1}}{1}}^2;$

\noindent
In particular, in both cases,
$d= 4$ and  the embedding is projectively normal.
\end{pr}

\noindent
{\em Proof.} Let   $F\simeq \pn{n-3}$ be any fiber of the scroll. We cut
(\ref{deg2dpf5}) with $F \cdot L^{n-5}$ and solve in $d$. We get $d=4$,
so that the  structure of $(X,L)$ is  given by
Theorem 8.10.1 of \ci{be-so-book}.

\noindent
In both cases it is easy to write down explicit equations for
the morphism associated with
$|\xi|$; we can check directly that $\xi$ is very ample, that the image
lies in a smooth quadric and that the embedding is projectively normal.
\blacksquare

\begin{rmk}
\label{alsoscrolloverthreefold}
{\rm Case  (\ref{pbundlesovercurves}.1) 
above is the Segre embedding of $\pn{1}\times\pn{3}$. It is a scroll 
over a curve
if we look at the first projection. If we look at the second projection
it is  a scroll over $\pn{3}$ with  associated vector bundle
${\odixl{\pn{3}}{1}}^{2}$.
Case (\ref{pbundlesovercurves}.2) is  a general hyperplane section of 
(\ref{pbundlesovercurves}.1);
the natural morphism onto $\pn{3}$ exhibits $X$ as the blow up of
$\pn{3}$ along a line.}
\end{rmk}

\subsection{Threefolds on $\Q{5}$ which are scrolls over surfaces} 
\label{scrollsoversurfaces}

\begin{lm}
\label{pbundlesoversurfaces}
Let $X\subseteq \Q{5}$ be  a codimension two  scroll over a surface $Y$.
Then  either $d=8$ or  we have:

\noindent
$g-1=\frac{1}{8}d(d-6)$,

\noindent
$\chi (\odix{X})=\chi (\odix{S})=\frac{1}{144}(d^3-18d^2+96d)$,

\noindent
$e_1^2=\frac{3}{2}d$, \quad  $e_2=\frac{d}{2}$,

\noindent
$K_Y\sim_n \frac{1}{6}(d-12) e_1$.
\end{lm}

\noindent
{\em Proof.} 
 We follow closely a  procedure which can be found  in \ci{ottp5}. 
By (\ref{euler}) we get:

\noindent
 $x_1=2L -e_1+y_1$;

\noindent
 $x_2= 2Ly_1 -e_1y_1+y_2$;

\noindent
 $x_3=2y_2$.
 
\noindent
We plug the above equalities in (\ref{deg2dpf5}) and (\ref{deg3dpf5}) 
and get
the following two equations:
\begin{equation}
\label{A}
(5-\frac{d}{2})L^2 + L\cdot e_1 -3L\cdot y_1 +
 e_1^2 + y_1^2 -e_1\cdot y_1 -y_2=0;
\end{equation}
\begin{equation}
\label{B}
\frac{d^2}{2}-2d + (d-8)L^2\cdot e_1 - (d-8)L^2\cdot y_1 -4L\cdot e_1^2 -
4L\cdot y_1^2 -2L\cdot y_2 +8L\cdot e_1 \cdot y_1 =0.
\end{equation}

\noindent
We cut (\ref{A}) and the tautological relation with
$L,$ $e_1$ and $y_1$ respectively. This way we get six 
elations which together
 with (\ref{B}) and the relation $L^3=d$ give a system of 
eight linear equations in the variables:
$v:=(L^3;\  L^2e_1;$ $ L^2y_1;$ $ Le_1^2;$ $ Le_1y_1;$ $ Ly_1^2;$ 
$ Le_2;$ $ Ly_2)$.
The matrix associated with the linear system is:

\[M:=  \left( \begin{array}{rrrrrrrr}
5-d/2 & 1 & -3 & 1 & -1 & 1 & 0 & -1  \\
0 & 5-d/2 & 0 & 1 & -3 & 0 & 0 & 0   \\
0 & 0 & 5-d/2 & 0 & 1 & -3 & 0 & 0  \\
d/2-2 & d-8 & -d+8 & -4 & 8 & -4 & 0 & -2 \\
1 & -1 & 0 & 0 & 0 & 0 & 1 & 0 \\
0 & 1 & 0 & -1 & 0 & 0 & 0 & 0 \\
0 & 0 & 1 & 0 & -1 & 0 & 0 & 0 \\
1 & 0 & 0 & 0 & 0 & 0 & 0 & 0
\end{array} 
\right) \]
and the linear system can be expressed as $Mv^t=(0,0,0,0,0,0,0,d)$.

\noindent
Since ${\rm det}\, M=72-9d$, the above  system of equations  
has a unique solution
if and only if
$d\neq 8$.

\noindent
Let us assume
$d\neq 8$.
Then the unique solution is:
\begin{eqnarray}
\label{solutiondnot8}
& &\{L^3;\ L^2e_1 ;\ L^2y_1 ;\  Le_1^2  ;\   Le_1y_1; 
\ Ly_1^2 ;\  Le_2;\   Ly_2 \} =   \\
 &  &\{  d;\ \frac{3d}{2} ;\  \frac{d}{4}(12-d);\  \frac{3d}{2}; 
\  \frac{d}{4}(12-d);\  \frac{d^3}{24}-d^2 +6d;\  \frac{d}{2};
\  \frac{d^3}{24}-\frac{d^2}{2}+2d\}. \nonumber
\end{eqnarray} 

\noindent
 We can use (\ref{solutiondnot8}) to
compute the genus of a general curve section, C, of $X$. 
This genus equals the arithmetic genus of
the line bundle $e_1$ on $Y$; we get
$$
2(g-1)= -e_1 y_1 + e_1^2=-Le_1y_1 + Le_1^2=\frac{1}{4}d(d-6).
$$
An analogous computation gives
\begin{equation}
\label{chitutti}
\chi (\odix{X})=\chi (\odix{Y}) = (1/144)(d^3-18d^2+96d)= \chi
(\odix{S}),
\end{equation}
where the first equality
is a standard fact about projective bundles
which can be proved using the Leray Spectral sequence 
and the last one follows
from the fact that $S$ is birationally equivalent to $Y$.

\noindent
To prove that $K_Y$ is numerically equivalent
to a rational multiple of $e_1$, we use Hodge index theorem for 
the surface
$Y$: by (\ref{solutiondnot8}), 
$K_Y^2 e_1^2=  (K_Y\cdot e_1)^2$, so that
$K_Y \sim_n q\,e_1$ for some rational number $q$ which is 
straightforward to compute.
\blacksquare

\begin{lm}
\label{d=6...42}
Let $X$ be a threefold scroll over a surface on $\Q{5}$.
Then  $d\leq 42$.
Moreover, if a general curve section, $C$, is contained in another
quadric hypersurface of \,$\pn{6}$, then $d\leq 12$. 
\end{lm}

\noindent
{\em Proof.}
By  Lemma \ref{pbundlesoversurfaces}
we have $g-1=(1/8)d(d-6)$.

\noindent
Assume that  a general curve section $C$ is not contained in any surface,
in the corresponding $\Q{3}$, of degree strictly less than
 $2\cdot 7$.
Then Proposition \ref{boundasep} implies $d\leq 42$. 

\noindent
Assume $C$ is contained in a surface ${\cal S}\subseteq \Q{3}$ of 
degree $2\cdot 6$. 
Proposition \ref{coreasybound} implies $d\leq 27$. 
\noindent
The same argument repeated for surfaces of degrees
$2\cdot 5$,
$2\cdot 4$ and $2\cdot 3$ gives  $d\leq 18$ in all three cases.

\noindent
Let us  assume that $C$ is contained in a surface of degree
$2\cdot 2$ and that $d>8$;    by Proposition \ref{roth},
 $X$ is contained in another quadric hypersurface of
$\pn{6}$. We now prove that, under the above  assumptions
on $C$, $d\leq 12$.
We plug $\s =2$ and the values of 
$\chi(\odix{S})$ and $g-1$, from Lemma \ref{pbundlesoversurfaces}, 
in inequality
(\ref{55}); we get
\begin{equation}
\label{pbundles=255}
-\frac{1}{288}d(d+6)(d-12)\geq -\frac{1}{12}[ \frac{1}{2}dL^2 -
 (K_X + 5L)^2]D.
\end{equation}
By Lemma \ref{pbundlesoversurfaces}  we have
$$
K_X=-2L +\frac{1}{6}(d-6)e_1.
$$
We now plug the above expression for $K_X$ in (\ref{pbundles=255})
using  the following relations
$e_2=(d/2)f$, $L^2D-Le_1D+e_2D=0$, where $f$ is a fiber
of the scroll. After simplifications the result is
$$
-d(d+6)(d-12)\geq 12(d+6)Le_1D +
d(d+6)(d-12)Df.
$$
Since $Le_1D\geq 0$ and $Df\geq 0$ we get $d\leq 12$.
Moreover, if $d=12$ then $D$ must be empty.

\noindent
Finally if $C$ were contained in a surface of degree $2\cdot 1$
then the same would be true for $X$,  by Theorem \ref{roth}.
But then $X$ would be a scroll on a quadric $\hat\Q{4}$ of $\pn{5}$ with 
at most one singular point.  Weil and Cartier
divisors coincide on $\hat\Q{4}$ and $Pic\,(\hat\Q{4})\simeq \zed$
by \ci{ha} II.6 Ex. 6.5. It would follow that $X$ is a complete 
intersection, a contradiction.  
\blacksquare

\begin{lm} 
\label{d=6812}
Let $X$ be a threefold scroll over a surface on $\Q{5}$. Then 
$d=6,$ $8$ or $12$.
\end{lm}

\medskip
\noindent
{\em Proof.} By Lemma \ref{d=6...42}, $d\leq 42$; by Lemma
\ref{pbundlesoversurfaces}, since the invariants there given must 
be integers
we see that the only possibilities for
the  pairs $(d,g)$ with $d>12$ are 
$(18,28)$, $(24,55)$, $(30,91)$, $(36,136)$ and $(42,190)$.
We prove that the cases $d=18,$ $24,$ $30,$ $36,$ $42$
cannot occur.

\noindent
Let $C\subseteq \Q{3} \subseteq \pn{4}$ be the general curve section of
$X$ and $\G \subseteq \Q{2} \subseteq \pn{3}$ be the general 
hyperplane section
of $C$. We denote by $h_C(i):= h^0( \odixl{\pn{4}}{i} ) -
h^0( {\cal I}_{C,\pn{4}}(i))$
the Hilbert function of $C\subseteq \pn{4}$ and by
$h_{\G}(i):= h^0( \odixl{\pn{3}}{i} ) -
 h^0( {\cal I}_{\G,\pn{3}} (i))$
the Hilbert function of $\G \subseteq \pn{3}$. 
Clearly $h_C(i)\leq h^0(\odixl{C}{i})$, 
for every $i$.

\noindent
{\em The case $(18,28)$.}  By Riemann-Roch
and Serre Duality we have $h^0(\odixl{C}{i})= 18i -27$
for $i\geq 4$; in particular
$h^0(\odixl{C}{4})= 45$.
 $C$ cannot be contained in another quadric of $\pn{4}$,
since otherwise, by Proposition \ref{roth} and Lemma \ref{d=6...42},
$d\leq 12$.
$C$ cannot be contained in an integral cubic of $\pn{4}$, otherwise,
we would get that the genus would be maximal with respect
to the bound prescribed by Proposition \ref{1.4}  and, since 
$\epsilon=0$, 
$C$
would be a complete intersection, forcing $X$ to be one too 
in view of  \ci{ha}, III.9 Ex. 9.6.
For the same reason   $C$ cannot be contained in an integral quartic of
 $\pn{4}$. 
It follows that there are no quartic hypersurfaces containing
$C$ except for the ones which are the union
of $\Q{3}$ with another quadric; in particular
$h_C(4)=55$. We get $55=h_C(4)\leq h^0(\odixl{C}{4})=45$, a contradiction.
The case $d=18$ cannot occur.

\smallskip
\noindent
{\em The case} $(24,55)$. As in the previous case we deduce that
$C$ is contained in a unique quadric of $\pn{4}$, $C$ is not contained in
 any integral
cubic or quartic of $\pn{4}$. This gives $h_C(4)=55$. As before
$h^0(\odixl{C}{5})=66$. By \ci{harrismont} Lemma 3.1 we have
$h_C(5) \geq h_C(4) + h_{\G}(5)$ and by \ci{harrismont}  Lemma 3.4
we also have that $h_{\G}(5)\geq 16$.  It follows that
$ 55 + 16 \leq h_C(5)  \leq  h^0(\odixl{C}{5})=66$, a contradiction.
The case $d=24$ cannot occur.

\smallskip
\noindent
{\em The cases $(30,91)$, $(36,136)$ and $(42,190)$}.
They are treated as the case $d=24$. In the first case
$h^0(\odixl{C}{7})=120$ and the only hypersurfaces of degree seven of 
$\pn{4}$
which contain $C$, contain $\Q{3}$, so that
$h_C(7)=140$, again a contradiction.
In the second case 
$h^0(\odixl{C}{8})=153$ and the only hypersurfaces of degree eight of 
$\pn{4}$
which contain $C$, contain $\Q{3}$, so that
$h_C(8)=289$, again a contradiction.
In the last case 
$h^0(\odixl{C}{10})=231$ and the only hypersurfaces of degree nine of 
$\pn{4}$
which contain $C$, contain $\Q{3}$, so that
$h_C(9)=385$. In particular $h_C(10)>385$, by  \ci{harrismont}
Lemma 3.1,  again a contradiction.
\blacksquare

\bigskip

The proof of the following is independent of
Theorem \ref{classificationd<12}.

\begin{pr}
\label{baseofscroll}
Let things be as in Lemma
{\rm \ref{d=6812}}.
If $d=6$ then $Y\simeq \pn{2}$; if $d=8$ then $Y\simeq \Q{2}$; if
$d=12$ then $Y$ is a minimal $K3$ surface.
\end{pr}

\noindent
{\em Proof.} Let $d=6$. By Lemma \ref{pbundlesoversurfaces}
$-K_Y$ is ample and $K_Y^2=9$; by the classification of Del Pezzo surfaces
we conclude that $Y\simeq \pn{2}$. Let $d=8$. The  proof
of Proposition \ref{nop2inq6} gives $Y \simeq \Q{2}$.
Let $d=12$. By Lemma \ref{pbundlesoversurfaces} $K_Y$ is numerically 
trivial, 
so that
$Y$ is a minimal model. (\ref{chitutti}) prescribes $\chi ({\cal 
O}_{Y})=2$, 
so that,
by the Enriques-Kodaira classification, $Y$ is a $K3$ surface.
\blacksquare

\subsubsection{The case of Type O)}
\label{d=12}

The purpose of this sections is twofold.
First we give an example of a scroll  of Type O), making the list of Theorem
\ref{maintm} effective.
Then we collect information about  the arbitrary  variety of this type.

Let $X$ be of Type  O), $\beta_{i,j}:=h^i( \Il{X}{\Q{5}}{j})$ and
  $\s_i:=h^i(\Il{X}{\Q{5}}{-1}\otimes {\cal S})$. The  sheaves $\Psi_i$, 
first introduced by
Kapranov, are defined
in  \ci{a-ott}.

\begin{tm}
\label{esempium}
Let
$
{\cal S}^7 \stackrel{\phi}\to {\Psi}_3 \oplus \odix{\Q{5}}^3 
$
be a  general  morphism. Then,  $\phi$ is injective,  $X:=D_{27}({\phi})$ 
is a variety of 
{\rm Type O)} and we have a resolution of the form
$$
0 \to {\cal S}^7 \stackrel{\phi}
\to {\Psi}_3 \oplus \odix{\Q{5}}^3 \to \Il{X}{\Q{5}}{3} \to 0.
$$
\end{tm}

We will prove this theorem after Proposition \ref{example}.
First we determine  some   properties of the arbitrary variety of Type O).
 
\begin{pr} 
\label{invd=12}
Let  $X\subseteq \Q{5}$ be of {\rm Type O)}.

\noindent
Then:
$$
g=10;\, \,  K_XL^2=-6; \,  \, K^2L=-6;\, \, K_X^3=12; \, \, c_2(X)L=24,
\, \, c_1({\cal E})^2=18, \, \, c_2({\cal E})=6.
$$

\smallskip
\noindent
The cohomology of $\odixl{X}{l}$:
$$
h^1({\odixl{X}{t}})=0, \,  \forall t \in {\zed};
$$
$$
h^2({\odixl{X}{t}})=0, \, \forall t \in {\zed} \, \,  except \, \, 
for \, \, h^2(\odix{X})=1; 
$$
$$
h^3({\odixl{X}{t}})=0, \, \forall t\geq -1.
$$
The following is the Beilinson-Kapranov $E^{p,q}_1$ table
for the sheaf $\Il{X}{\Q{5}}{3}$; see {\rm \ci{a-ott}} {\rm Theorem
5.6}. A letter {\rm a} on the
left of a  vector bundle $B$ means $B$ direct sum with itself {\rm a} times.
$$ 
\begin{array}{cccccc}
0 & 0 & 0 & 0 & 0 & 0    \\
7{\cal S} & 0 & 0 & 0 & 0 & 0     \\
0  & 0 & \Psi_3 & 0 & 0 & 0    \\
0 & 0 & 0 & 0  & 0   & 0    \\
0 & 0 & 0 & 0  & \beta_{1,2} \Psi_1  & \beta_{1,3} \odix{X}   \\
0 & 0 & 0 & 0 & \beta_{0,2} \Psi_1 & \beta_{0,3} \odix{X},   
\end{array}
$$ 
where
$\beta_{0,2} - \beta_{1,2} =0$, 
 $\beta_{0,3} - \beta_{1,3}=3$. 

\noindent
 Either
$\beta_{2,0}=1$ and $\beta_{0,3}=7$ or $\beta_{2,0}=0$
and $ 0 \leq \beta_{1,3} \leq 21$.

\noindent
If $\beta_{2,0}=1$ then  $\Il{X}{\Q{5}}{3}$ can be expressed
as the cohomology of a  monad of the form {\rm (see \ci{o-s-s} 
for the definition
of monads)}:
\begin{equation}
\label{II2}
0 \to 7{\cal S} \stackrel{m_1}{\to} \Psi_3 
\oplus \Psi_1 \oplus \odixl{\Q{5}}{1} 
\stackrel{n_1}{\to} 4 \odix{\Q{5}} 
\to 0
\end{equation}
If $\beta_{2,0}=0$ then $\Il{X}{\Q{5}}{3}$ can be expressed
as the cohomology of a  monad of the form:
\begin{equation}
\label{I3}
0 \to 7{\cal S} \stackrel{m_2}{\to} \Psi_3  
\oplus 3 \odix{\Q{5}}\oplus \beta_{1,3}
\odix{\Q{5}} \stackrel{n_2}{\to} \beta_{1,3} \odix{\Q{5}} \to 0.
\end{equation}

\end{pr}

\noindent
{\em Proof.}
The first list of invariants can be read off from Lemma   
ref{scrollsoversurfaces}
when $d=12$.

\noindent
As to $h^2(\odixl{X}{t})$ we argue as follows. 
Via the projection formula and Leray Spectral Sequence,
$h^2(\odixl{X}{t})=0$, $\forall t<0$ and for the same reason
$h^2(\odix{X})=1$. Since $K_Y$ is trivial,
Leray Spectral Sequence and
Le Potier's Vanishing Theorem \ci{sh-so} give
 $h^i(\odixl{X}{1} )=h^i(\odixl{Y}{\cal E})=
h^i(\odixl{Y}{K_Y \otimes {\cal E}}) =0$,
$\forall $ $i\geq 2$. 
By Serre Duality and  the fact that $L_{|S}(K_S -mL_{|S})=6 - 12m$
 we see that
$h^2(\odixl{S}{m}) =0$, $\forall m\geq  1$;  we conclude 
for $h^2$ by an easy induction using
the sequences 
\begin{equation}
\label{hyperseq}
0 \to \odixl{X}{m-1} \to \odixl{X}{m} \to \odixl{S}{m} \to 0.
\end{equation}
The vanishings of the $h^3$'s are obvious consequences of
Serre Duality.

\noindent
$h^1(\odixl{X}{t}=0, \, \forall t<0$ by Kodaira vanishing. 
For $t=0$ the vanishing follows from
$h^1(\odix{Y})=0$. Since $X\subseteq \pn{6}$
is linearly normal by  a result of Fujita's (cf. \ci{okcod3} , \S 4) and
$\chi(\odixl{X}{1})=7$ by Riemann-Roch, 
we have $h^1(\odixl{X}{1})=0$. To prove the remaining vanishings
for $h^1(\odixl{X}{t})$ we argue by induction using the long cohomology 
sequences associated
with the sequences
(\ref{hyperseq}), the analogue ones obtained by replacing $X$ 
and $S$ by $S$ and $C$ 
(a general curve section of $S$) and observing that the linear systems
$|\odixl{C}{t}|$ are non-special for $t\geq 2$.

\noindent
The Beilinson-Kapranov table is obtained as follows.
The vanishings $\beta_{i,j}=0$ for $i=2,\, 3,\, 4,\,5$ and 
$j=-1,\,0,\,   1,\,
2,\, 3$, except for $\beta_{3,0}=h^0(\odix{X})=1$, are obtained
by taking the cohomology of  the exact sequences 
\begin{equation}
\label{IOX}
0
\rightarrow
\Il{X}{\Q{5}}{l}
\rightarrow
\odixl{\Q{5}}{l}
\rightarrow
\odixl{X}{l}
\rightarrow
0
\end{equation}
and plugging the above values for the cohomology of $\odixl{X}{t}$.
For the same reason $\beta_{i,j}=0$ for $i=0, \, 1, \, 2$, $j=-1,\, 0$.
$\beta_{0,1}$ is zero because $X$ cannot be degenerate (see the proof
of Lemma \ref{d=6...42}).
$\beta_{1,1}=0$ since $X\subseteq \pn{6}$ is linearly normal.
The relations on the remaining $\beta$'s come from the shape of 
the Hilbert polynomial
$$
\chi (\Il{X}{\Q{5}}{t})={\frac {1}{60}}t^5+{\frac {5}{24}}t^{4}
-t^{3}+{\frac {19}{24
}}t^{2}+{\frac {59}{60}}t-1,
$$
which vanishes for $t=-1,\, 1,\, 2,$ and has value three for $t=3$.

\noindent
Because of how  this spectral sequence works ($E_{\infty}=E_6$
and
$E^{p,q}_{\infty}\simeq \{0\}$ for $p+q\not=0$), 
 we see that
$\s_i=0$, for $i=0,\,1,\,2$.

\noindent
${\s}_5=0$ by observing the cohomology  of (\ref{IOX}) 
tensored with ${\cal S}$.
We use the same sequences, together with Riemann-Roch for
$\cal S$ and for ${\cal S}_{|X}$ to get
$$
\chi ({\cal S}(t))=\frac{1}{15}t(t+1)(t+2)(t+3)(t+4),
$$
and
$$
\chi ({\cal S}_{|X} (t))=8t^3 -6t^2 +7;
$$
 it follows that
$
-\s_3 + \s_4=\chi(\I{X}{\Q{5}}\otimes {\cal S}(-1))=7.
$

\noindent
We now prove that $\s_3 =0$. There is at most
one nontrivial differential from $\s_3{\cal S}$, namely 
the one that hits
$E^{-2,0}_4$. On the other hand $E^{-2,0}_4=E^{-2,0}_2=\mbox{Ker}d^{-2,0}_1$.
It is enough to show that the last group is trivial.
We consider two cases. The former is when $\beta_{0,2}=0$; in this case
$\s_3$ is clearly zero. The latter is when $\beta_{0,2}\not=0$. 
Then $\beta_{0,2}=1$
otherwise $X$ would have $d\leq 8$, the degree of the intersection 
on $\Q{5}$
of two hypersurfaces of degree four. $\beta_{0,3}=7$ otherwise
$X$ would be a complete intersection on $\pn{6}$ of type $(2,2,3)$,
 a contradiction.
By Kapranov's  explicit resolution of the diagonal on $\Q{n}\times \Q{n}$,
 see
\ci {a-ott},  
  we infer that $d^{-2,0}_1$  coincides with the injection
$(\Psi_1\simeq)\, {\Omega^1(1)_{\pn{6}}}_{|\Q{5}} \to \odix{\Q{5}}^7$ 
obtained by restricting
the Euler sequence on $\pn{6}$ to $\Q{5}$. It follows that $\s_3=0$ 
and therefore $\s_4 =7$.

\noindent
If $\beta_{2,0}\not=0$, we have seen above that $\beta_{2,0}=1$, 
$\beta_{3,0}=7$
and that $d^{-2,0}_1$ coincides with the injection
$(\Psi_1\simeq)\, {\Omega^1(1)_{\pn{6}}}_{|\Q{5}} \to \odix{\Q{5}}^7$ 
whose cokernel is
$\odixl{\Q{5}}{1}$. The statement associated with (\ref{II2}) follows 
from 
\ci{a-ott}.

\noindent
Similarly, we see that the statement associated with (\ref{I3}) holds
when $\beta_{2,0}=0$.
Since the morphism 
$n_2$ is  trivial on $(3+\beta_{1,3})\odix{\Q{5}}$, the restriction  
$\nu:={n_2}_{|\Psi_3}$
is surjective. Recall that the rank of
$\Psi_3$ is $26$. If $\beta_{1,3}>21$, then the
kernel of the map $\nu$ would be a locally free sheaf
of rank $r<5$ with fifth Chern class $c_5=c_5(\Psi_3)\not=0$, a 
contradiction.
\blacksquare

\bigskip

The following is  essentially due to Peskine and Szpiro; see \ci{okcod3},
 \S1.

\begin{lm}
\label{pesk}
Let $X$ be a codimension two nonsingular subvariety of a nonsingular 
variety $Z$ of
dimension $n\leq 5$ and $L_i$, $i=1,\,2$, two line bundles on $Z$ such that
the sheaves
$\Il{X}{Z}{L_i}$ are globally generated on $Z$.
Let $s_i \in $ $H^0(\Il{X}{Z}{L_i})$ be the choice of two general sections
and $V_i$ the two effective divisors associated with the $s_i$.
Then $V_1 \cap V_2= X \cup Y,$ as schemes,  where $Y$ is nonsingular. 
\end{lm}

We now give a family of  examples of degree $d=12$ scrolls on $\Q{5}$. 

\begin{pr}
\label{example}
Let  $\rho: {\cal C}^{\vee} \to \odixl{\Q{5}}{2}^3$ be a generic morphism.
Then $X:=D_1(\rho)$ is a variety of {\rm Type O)} such that
$\Il{X}{\Q{5}}{3}$ is generated by global sections; $X$ is linked
to a variety $X'$ of \,{\rm Type F)} via  the complete intersection
of two general elements of $|{\Il{X}{\Q{5}}{3}}|$.

\noindent 
Conversely, if $X\subseteq \Q{5}$ is of {\rm Type O)}
and $\Il{X}{\Q{5}}{3}$ is generated by  global sections,
then $X=D_1(\rho)$ for some $\rho$ as above and $X$ is linked as above to 
a variety of {\rm Type F)}.
\end{pr}

\noindent
{\em Proof.} For facts about Cayley bundles,  see
 Fact \ref{cayley}. ${\cal C}(2)^3$ is generated by global sections and
 Fact \ref{hartshorne} implies that  $X:=D_1(\rho)$ is a codimension two 
nonsingular 
subvariety of
$\Q{5}$ and that we have  the following exact sequence
\begin{equation}
\label{alceste}
 0 \to {\cal C}^{\vee}(-2) \stackrel{\rho}{\to} \odix{\Q{5}}^3 
\to \Il{X}{\Q{5}}{3} \to 0.
\end{equation}
We compute the total Chern class of  \,$\Il{X}{\Q{5}}{3}$ via 
(\ref{alceste}): 
$ 1+3\,h+6\,h^{2}+9\,h^{3}+9\,h^{4}$.

\noindent
We compare it with  Lemma
\ref{lmccis}\,:

\noindent
$
\gamma_1=3,$\,  $\gamma_2=6=\frac{1}{2}d,$\, 
$\gamma_3=\frac{1}{2}(K_X+2L)L^2=9,$\,
$ \gamma_4=
\frac{1}{2}(K_X +2L)^2L=9,$\, $ \gamma_5=\frac{1}{2}(K_X+2L)^3=0,$

\noindent
where $L$ denotes $\odixl{\pn{6}}{1}_{|X}$.
It follows that $X$ has degree $d=12$. 
By \ci{be-so-book} Proposition 7.2.2, $K_X+2L$ 
is generated by global sections
since $(X,L)$ cannot be  isomorphic to either
$(\pn{3}, \odixl{\pn{3}}{1})$, $(\Q{3}, \odixl{\Q{3}}{1})$ or to 
a scroll
over a curve since they all have degree $d=4$ by Proposition
\ref{pbundlesovercurves}. The fact that $\gamma_5=0$ implies 
  that $K_X+2L$ cannot  be big and the fact that $\gamma_4\not= 0$
implies that the Kodaira dimension $\kappa (K_X+2L)=2$, so that, by
\ci{be-so-book} Theorem 7.3.2, 
$(X,L)$ is an adjunction theoretic scroll over a surface and,
by
\ci{be-so-book}  Proposition 14.1.3, it is actually a scroll 
in our sense. By
Lemma \ref{baseofscroll}, $X$ is a degree $d=12$ scroll over
a $K3$ surface. The linking part is proved using Lemma
\ref{pesk} to produce a $X'$ of  degree $d'=6$ and by observing that
the mapping cone construction yields  a  resolution
for $\Il{X'}{\Q{5}}{3}$ which coincides with the one of a variety of 
Type F).

\noindent
The converse is proved in a similar way.
\blacksquare

\bigskip
\noindent
 {\bf Proof of Theorem}  {\bf \ref{esempium}}.
For the varieties constructed in Proposition \ref{example} we have,
by (\ref{alceste}), that $\beta_{2,0}=0$ and $\beta_{0,3}=3$. 
By looking at the display
of the monad (\ref{I3}) with these invariants we get
 the desired resolution for the ideal sheaves of these varieties. 
It also follows that,
 for the generic morphism $\phi$,
$\phi$  is injective and $X:=D_{27}({\phi})$ is of  Type O)
 as in the proof of
Proposition \ref{example}.
\blacksquare

\subsection{4-folds which are scrolls on $\Q{6}$}
\label{scrollsinq6}

\begin{pr}
\label{nop2inq6}
There are no fourfolds scrolls over surfaces on $\Q{6}$.
\end{pr}

\noindent
{\em Proof.}
By contradiction, assume that $X^4$ is such a scroll.
Cutting
(\ref{deg2dpf5}) with a fiber $F\simeq \pn{2}$, we get $d=8$. 
We take a general
hyperplane section and  obtain a scroll, $X$,
 on $\Q{5}$ so that the previous analysis applies.  In fact,
if a special fiber, $F$, were isomorphic to $\pn{2}$, then $F_{|F}\simeq
\odixl{\pn{2}}{-1}$ so that we would have a contraction morphism
$\eta: X\to X'$ and the structural morphism $p:X \to Y$  would factor 
through
$\eta$ violating 
 the upper semicontinuity of the dimension of the fibers.

\noindent
We solve the linear system contained in the proof of 
Lemma \ref{scrollsoversurfaces} for $d=8$ and  we get that the solutions
depend on one additional  parameter $t$:
\begin{eqnarray}
\label{solutiond8}
& &\{L^3;L^2e_1 ;L^2y_1 ; Le_1^2  ;  Le_1y_1; 
Ly_1^2 ; Le_2;  Ly_2 \} =  \\
 &  &\{  8; 36-\frac{9}{2}t ;24-3t;36-\frac{9}{2}t; 
 24-3t; 16-2t;28- \frac{9}{2}t;
 t \}. \nonumber
\end{eqnarray} 
We observe that
$K_Y^2=(4/9)e_1^2$ and that  $K_Y\cdot e_1=-(2/3)e_1^2$. Since
$e_1$ is ample, the Hodge Index Theorem implies that 
 $3K_Y\sim_n -2e_1$. It follows that $Y$ has to be a Del Pezzo surface.
On such a $Y$,
 numerical and rational equivalence
coincide and $3K_Y$ is not divisible by 2 unless $Y$ is a smooth
$\Q{2}$. In this case $t=b_2=4$.
In particular $\deg e_2 = 10$ and $g=4$.  
 By \ci{io3} this case cannot occur
if $\dim X\geq 4$. 
This contradicts the existence of  scrolls over surfaces
in $\Q{6}$ of dimension four.
\blacksquare

\begin{pr}
\label{p1bundlesinq6}
The only  scroll over   a threefold on $\Q{6}$ is \,
$\pn{1}\times \pn{3}$ embedded with the Segre embedding.

\end{pr}

\noindent
{\em Proof.} The  proof runs along the lines of
Lemma  \ref{pbundlesoversurfaces}.
Using (\ref{euler}) we compute the   Chern classes of $X$:

\noindent
$x_1=2L -e_1 +y_1$;

\noindent
$x_2= 2Ly_1 -e_1y_1 +y_2$;

\noindent
$x_3= 2Ly_2 - e_1y_2 +y_3$;

\noindent
$x_4=2Ly_3$.

\noindent
After having plugged these relations in (\ref{deg2dpf5}), 
(\ref{deg3dpf5}) and 
(\ref{deg4dpf})
we get, respectively:
\begin{equation}
\label{n2}
( \frac{1}{2} d -8)L^2 + L(4y_1-2e_1)-e_1^2+e_1y_1-y_1^2+y_2=0,
\end{equation}
\begin{eqnarray}
\label{n3}
& & 8L^3+ L^2(4e_1-8y_1)+L(-2e_1y_1+4y_1^2-4y_2)+  \\
& & e_1^3-e_1^2y_1+e_1y_1^2-e_1y_2-y_1^3
+2y_1y_2-y_3=0,  \nonumber
\end{eqnarray}
\begin{eqnarray}
\label{n4}
& & 6L^4-8L^3y_1+L^2(4e_1^2-4e_1y_1+8y_1^2-8y_2) +    \\
& & L(-2e_1^3+2e_1y_1^2-2e_1y_2-4y_1^3
+8y_1y_2-4y_3)=0. \nonumber
\end{eqnarray}
We cut the tautological relation and (\ref{n2}) with the following classes:
$L^2$, $Le_1$, $Ly_1$, $e_1^2$, $ y_1^2$, $e_1y_1$, $e_2$ and $y_2$; we cut
(\ref{n3}) with $L$, $e_1$ and $y_1$. Considering also (\ref{n4}) we get 
a total of twenty linear equations in the seventeen variables:
$L^4$, $L^3e_1$, $L^3y_1$, $L^2e_1^2$, $L^2e_1y_1$, $L^2e_2$, $L^2y_1^2$,
$L^2y_2$, $ Le_1^3$, $Le_1^2y_1$, $Le_1e_2$, $Le_1y_1^2$,
$Le_1y_2$, $Le_2y_1$, $Ly_1^3$,
$Ly_1y_2$ and $Ly_3$.
We leave out, on purpose, the condition $L^4=d$. We leave to 
the reader to check
that the resulting linear system  has a nontrivial solution
only for $d=4$.
If $d=4$ we use Theorem \ref{classificationd<12} to conclude.
\blacksquare

\bigskip
1991 {\em Mathematics Subject Classification}. 14C05, 14F05, 14H45, 
14J10, 14J28,
14J40, 14J60, 14M06,14M07, 14M17.

{\em Key words and phrases}. Classification; 
 liaison; low codimension; low degree; quadric; scroll; vector bundle.

\bigskip
AUTHOR'S ADDRESS

\medskip
\noindent
Mark Andrea de Cataldo,
Department of Mathematics,
Washington University in St. Louis,
Campus Box 1146,
St. Louis, Missouri 63130-4899.

\noindent
e-mail: mde@math.wustl.edu

\end{document}